\documentclass[letterpaper, 10 pt, conference]{ieeeconf}

\IEEEoverridecommandlockouts    
\usepackage{psfrag}
\usepackage{amsmath,amssymb,amsfonts,bbm}
\usepackage{latexsym}
\usepackage{graphicx}
\usepackage{colortbl}
\usepackage{fancyhdr}
\usepackage{epstopdf}
\usepackage{float}
\usepackage{hyperref}
\usepackage{booktabs}
\usepackage{enumerate}
\usepackage[ruled]{algorithm2e}
\usepackage{balance}
\usepackage{mathrsfs}
\usepackage{tikz}
\usepackage{pgfplots}
\pgfplotsset{compat=1.10}
\usepgfplotslibrary{fillbetween}

\usepackage{mathtools}

\newtheorem{theorem}{Theorem}
\newtheorem{definition}{Definition}
\newtheorem{proposition}{Proposition}
\newtheorem{lemma}{Lemma}
\newtheorem{corollary}{Corollary}

\newtheorem{assumption}{Assumption}

\newcommand{\R}{\mathbb{R}}

\newcommand{\D}{\mathbb{D}}

\newcommand{\N}{\mathbb{N}}

\newcommand{\mc}{\mathcal}

\newcommand{\bbS}{\mathbb{S}}

\newcommand{\cT}{\mathcal{T}}
\newcommand{\cA}{\mathcal{A}}

\newcommand{\fix}{\mathrm{fix}}
\newcommand{\zer}{\mathrm{zer}}
\newcommand{\Id}{\mathrm{Id}}

\newcommand{\bdr}{\mathrm{bdr}}
\newcommand{\bs}{\boldsymbol}

\newcommand\oprocendsymbol{\hbox{$\square$}}
\newcommand\oprocend{\relax\ifmmode\else\unskip\hfill\fi\oprocendsymbol}

\begin{document}

\title{On the convergence of discrete-time linear systems: A linear time-varying Mann iteration converges iff the operator is strictly pseudocontractive}

\author{
Giuseppe Belgioioso \and Filippo Fabiani \and Franco Blanchini \and Sergio Grammatico
\thanks{G. Belgioioso is with the Control Systems group, TU Eindhoven, The Netherlands. 
F. Fabiani is with the Department of Information Engineering, University of Pisa, Italy.
F. Blanchini is with the Department of Mathematics and Informatics, University of Udine, Italy.
S. Grammatico is with the Delft Center for Systems and Control (DCSC), TU Delft, The Netherlands.
E-mail addresses: \texttt{g.belgioioso@tue.nl}, \texttt{filippo.fabiani@ing.unipi.it}, \texttt{franco.blanchini@uniud.it}, \texttt{s.grammatico@tudelft.nl}. This work was partially supported by NWO under research projects OMEGA (grant n. 613.001.702) and P2P-TALES (grant n. 647.003.003).
\smallskip \newline
}
}

\thispagestyle{empty}
\pagestyle{empty}

\maketitle         

\begin{abstract}
We adopt an operator-theoretic perspective to study convergence of linear fixed-point iterations and discrete-time linear systems. We mainly focus on the so-called Krasnoselskij--Mann iteration 
$x(k + 1) = (1-\alpha_k) x(k) + \alpha_k A x(k)$, which is relevant for distributed computation in optimization and game theory, when $A$ is not available in a centralized way. 
We show that convergence to a vector in the kernel of $(I-A)$ is equivalent to strict pseudocontractiveness of the linear operator $x \mapsto A x$.
We also characterize some relevant operator-theoretic properties of linear operators via eigenvalue location and linear matrix inequalities. We apply the convergence conditions to multi-agent linear systems with vanishing step sizes, in particular, to linear consensus dynamics and equilibrium seeking in monotone linear-quadratic games.
\end{abstract}

\section{Introduction}

State convergence is the quintessential problem in multi-agent systems. In fact, multi-agent {consensus} and cooperation, distributed {optimization} and multi-player {game theory} revolve around the convergence of the state variables to an equilibrium, typically unknown a-priori. In distributed consensus problems, agents interact with their neighboring peers to collectively achieve global agreement on some value \cite{olfati-saber:fax:murray:10}. In distributed optimization, decision makers cooperate locally to agree on primal-dual variables that solve a global optimization problem \cite{nedic:ozdaglar:parrillo:10}. 
Similarly, in multi-player games, selfish decision makers exchange local or semi-global information to achieve an equilibrium for their inter-dependent optimization problems \cite{yi:pavel:17cdc}.
Applications of multi-agent systems with guaranteed convergence are indeed vast, e.g. include power systems \cite{dorfer:simpson-porco:bullo:16, dorfler:grammatico:17}, demand side management \cite{mohsenian-rad:10}, network congestion control \cite{jaina:walrand:10, barrera:garcia:15}, social networks \cite{Ghaderi2014, etesami:basar:15}, robotic and sensor networks \cite{martinez:bullo:cortes:frazzoli:07, sankovic:johansson:stipanovic:12}.

From a general mathematical perspective, the convergence problem is a fixed-point problem \cite{berinde}, or equivalently, a zero-finding problem \cite{bauschke:combettes}. For example, consensus in multi-agent systems is equivalent to finding a collective state in the kernel of the Laplacian matrix, i.e., in operator-theoretic terms, to finding a zero of the Laplacian, seen as a linear operator.

{Fixed-point theory} and {monotone operator theory} are then key to study convergence to multi-agent equilibria \cite{ryu:boyd:16}. For instance, Krasnoselskij--Mann fixed-point iterations have been adopted in aggregative game theory \cite{grammatico:parise:colombino:lygeros:16, grammatico:17}, monotone operator splitting methods in distributed convex optimization \cite{giselsson:boyd:17} and monotone game theory \cite{belgioioso:grammatico:17cdc,yi:pavel:17cdc,grammatico:17tcns}. The main feature of the available results is that \textit{sufficient} conditions on the problem data are typically proposed to ensure global convergence of fixed-point iterations applied on \textit{nonlinear} mappings, e.g. compositions of proximal or projection operators and linear averaging operators.

Differently from the literature, in this paper, we are interested in \textit{necessary and sufficient} conditions for convergence, hence we focus on the three most popular fixed-point iterations applied on \textit{linear} operators, that essentially are linear time-varying systems with special structure. The motivation is twofold. First, there are still several classes of multi-agent linear systems where convergence is the primary challenge, e.g. in distributed linear time-varying consensus dynamics with unknown graph connectivity (see Section \ref{sec:App-Cons}). Second, fixed-point iterations applied on linear operators can provide non-convergence certificates for multi-agent dynamics that arise from distributed convex optimization and monotone game theory (Section \ref{sec:App-Game}).

Our main contribution is to show that the  Krasnoselskij--Mann fixed-point iterations, possibly time-varying, applied on linear operators converge if and only if the associated matrix has certain spectral properties (Section \ref{sec:MainRes}). To motivate and achieve our main result, we adopt an operator-theoretic perspective and characterize some regularity properties of linear mappings via eigenvalue location and properties, and linear matrix inequalities (Section \ref{sec:OTPLM}). In Section \ref{sec:Concl},  we conclude the paper and indicate one future research direction.

\begin{figure}[!t]
\begin{center}
\includegraphics[width=1\columnwidth]{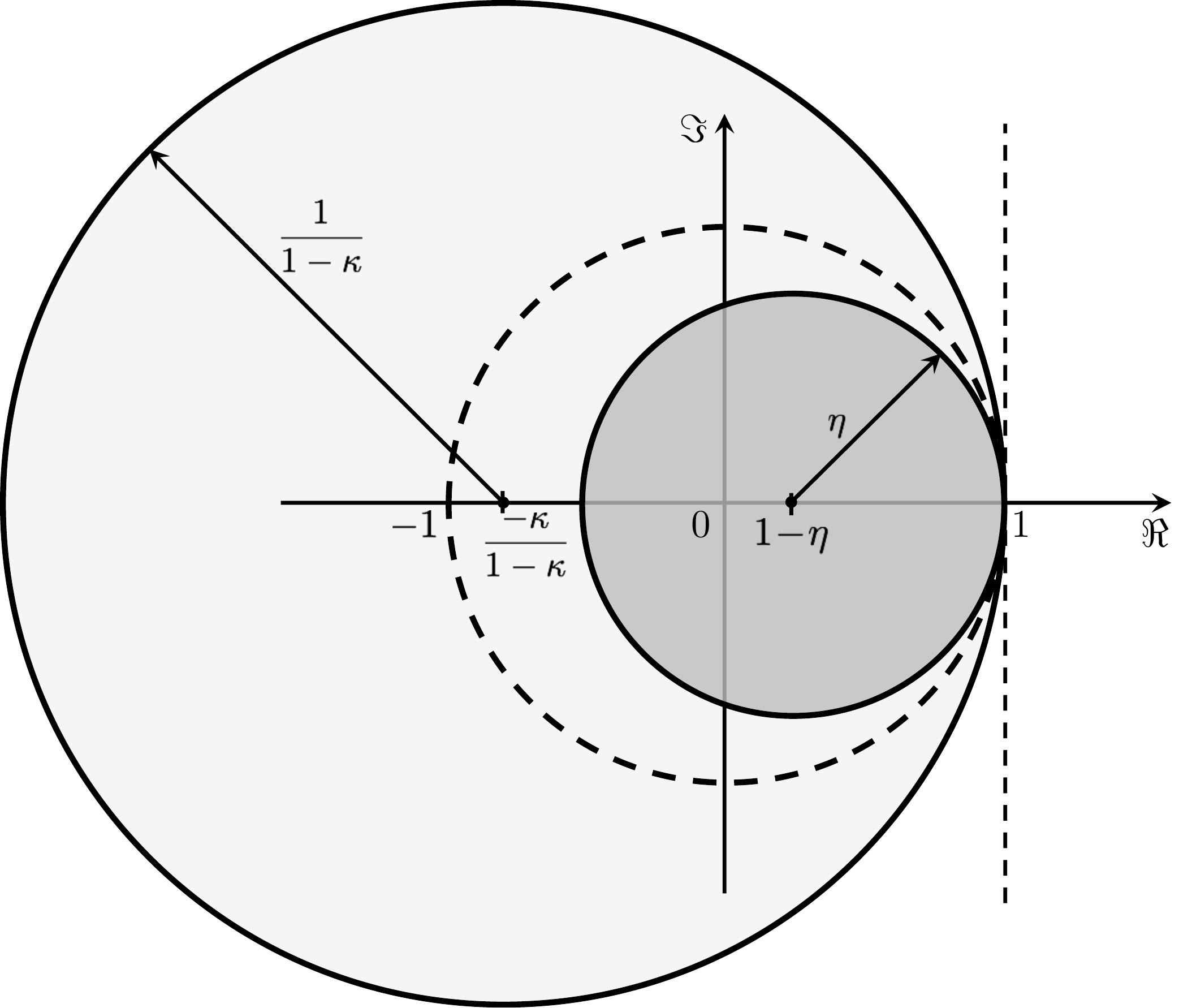}
\caption{Spectrum of a linear $\eta$-AVG mapping: Disk centered in $1-\eta$ with radius $\eta$, $\mathbb{D}_{\eta}$ (dark-grey disk). Spectrum of a linear $\kappa$-sPC mapping: Disk centered in $-\frac{\kappa }{1-\kappa }$ with radius $\frac{1}{1-\kappa }$, $\mathbb{D}_{\frac{1}{1-\kappa }}$ (light-grey disk).}
\label{fig:spc}
\end{center}
\end{figure}

\textit{Notation}: $\R$, $\R_{\geq 0}$ and $\mathbb{C}$ denote the set of real, non-negative real and complex numbers, respectively. $\D_r := \{ z \in \mathbb{C} \mid | z - (1-r) | \leq r \}$ denotes the disk of radius $r>0$ centered in $(1-r,0)$, see Fig. \ref{fig:spc} for some graphical examples. 
$\mathcal{H}\left( \left\| \cdot \right\| \right)$ denotes a finite-dimensional Hilbert space with norm $\left\| \cdot \right\|$. 
$\bbS^n_{\succ 0}$ is the set of positive definite symmetric matrices
and, for $P \in\bbS^n_{\succ 0}$, $\| x \|_P := \sqrt{ x^\top P x }$.
$\Id$ denotes the identity operator. 
$R(\cdot) := \left[ 
\begin{smallmatrix}
\cos(\cdot) & -\sin(\cdot) \\
\sin(\cdot) & \phantom{-}\cos(\cdot)
\end{smallmatrix}
\right]$ denotes the rotation operator. 
Given a mapping $T:\R^n \rightarrow \R^n$, $\textrm{fix}(T) := \left\{ x \in \R^n \mid x = T(x) \right\}$ denotes the set of fixed points, and $\textrm{zer}(T) := \left\{ x \in \R^n \mid 0 = T(x) \right\}$ the set of zeros. Given a matrix $A \in \R^{n \times n}$, $\textrm{ker}(A) := \left\{ x \in \R^n \mid 0 = A x \right\} = \textrm{zer}(A \, \cdot)$ denotes its kernel; $\Lambda(A)$ and $\rho(A)$ denote the spectrum and  the spectral radius of $A$, respectively. $\bs{0}_N$ and $\bs{1}_N$ denote vectors with $N$ elements all equal to $0$ and $1$, respectively.

%----------------------------------------------------------------
\section{Mathematical definitions} \label{sec:MB}

\subsection{Discrete-time linear systems}

In this paper, we consider discrete-time linear time-invariant systems,
\begin{equation}
\label{eq:system}
x(k+1) = A x(k)\,,
\end{equation}
and linear time-varying systems with special structure, i.e.,
\begin{equation}
\label{eq:tv-system}
x(k+1) = (1-\alpha_k) x(k) + \alpha_k A \, x(k)\,,
\end{equation}
for some positive sequence $\left( \alpha_k \right)_{k \in \N}$. Note that for $\alpha_k=1$ for all $k \in \N$, the system in \eqref{eq:tv-system} reduces to that in \eqref{eq:system}.

\subsection{System-theoretic definitions}

We are interested in the following notion of global convergence, i.e., convergence of the state solution, independently on the initial condition, to some vector.

\smallskip
\begin{definition}[Convergence] \label{def:convergence}
The system in \eqref{eq:tv-system} is \textit{convergent} if, for all $x(0) \in \R^n$, its solution $x(k)$ converges to some $\bar{x} \in \R^n$, i.e., $\displaystyle \lim_{k \rightarrow \infty} \left\| x(k) - \bar{x} \right\| = 0$.
\hfill $\square$
\end{definition}
\smallskip

Note that in Definition \ref{def:convergence}, the vector $\bar{x}$ can depend on the initial condition $x(0)$.
In the linear time-invariant case, \eqref{eq:system}, it is known that semi-convergence holds if and only if the eigenvalues of the $A$ matrix are strictly inside the unit disk and the eigenvalue in $1$, if present, must be semi-simple, as formalized next.
\smallskip
\begin{definition}[(Semi-) Simple eigenvalue]
An eigenvalue is \textit{semi-simple} if it has equal algebraic and geometric multiplicity. An eigenvalue is \textit{simple} if it has algebraic and geometric multiplicities both equal to $1$.
\hfill $\square$
\end{definition}
\smallskip

\smallskip
\begin{lemma} \label{lem:conv-eig}
The following statements are equivalent:
\begin{enumerate}
\item[i)] The system in \eqref{eq:system} is \textit{convergent};
\item[ii)] $\rho(A) \leq 1$ and the only eigenvalue on the unit disk is $1$, which is semi-simple.
\hfill $\square$
\end{enumerate}
\end{lemma}

%\begin{definition}[Weak Lyapunov function]
%A function $v: \R^n \rightarrow \R_{\geq 0}$ is a \textit{weak Lyapunov function} for the system in \eqref{eq:system} if $v( A x ) \leq v( x )$ for all $x \in \R^n$.
%\hfill $\square$
%\end{definition}

\subsection{Operator-theoretic definitions} \label{subsec:OTD}

With the aim to study convergence of the dynamics in \eqref{eq:system}, \eqref{eq:tv-system}, in this subsection, we introduce some key notions from operator theory in Hilbert spaces.

\begin{definition}[Lipschitz continuity]
\label{def:Lipschitz}
A mapping $\mc{T}: \R^n \rightarrow {\R}^n$ is $\ell$-\textit{Lipschitz} continuous in $\mathcal{H}\left( \left\| \cdot \right\| \right)$, with $\ell \geq 0$, if $\forall x,y \in \R^n$,  
$ \left\| \mc{T}(x) - \mc{T}(y)\right\| \leq \ell \left\| x-y\right\|.$
{\hfill $\square$}
\end{definition}
\smallskip

\begin{definition} \label{def:NE}
In $\mathcal{H}\left( \left\| \cdot \right\| \right)$, an $\ell$-Lipschitz continuous mapping $\mc{T}: \R^n \rightarrow {\R}^n$ is 
\begin{itemize}
\item $\ell$-\textit{Contractive} ($\ell$-\textup{CON}) if $\ell \in [0,1)$;
\smallskip
\item \textit{NonExpansive} (NE) if $\ell \in [0,1]$;
\smallskip
\item $\eta$-\textit{Averaged} ($\eta$-\textup{AVG}), with $\eta \in (0,1)$, if $\forall x,y \in \R^n$
\begin{multline} \label{eq:AVG-ineq}
 \left\| \mc{T}(x) - \mc{T}(y) \right\|^2  \leq \left\| x-y \right\|^2 \\
  - \tfrac{1-\eta}{\eta}\left\| \left( \textup{Id}-\mc{T} \right)(x) - \left(\textup{Id}-\mc{T} \right)(y) \right\|^2 \, ,
\end{multline}
or, equivalently, if there exists a nonexpansive mapping $\mathcal{B}: \R^n \rightarrow \R^n$ and $\eta \in (0,1)$ such that 
$$\mathcal{T} = (1-\eta) \textup{Id} + \eta \mathcal{B} \,.$$
\item $\kappa$-\textit{strictly Pseudo-Contractive} ($\kappa$-sPC), with $\kappa \in (0,1)$, if $\forall x,y \in \R^n$
\begin{multline} \label{eq:SPC-ineq}
\left\| \mc{T}(x) - \mc{T}(y) \right\|^2  \leq \left\| x-y \right\|^2 \\
  + \kappa \left\| \left( \textup{Id}-\mc{T} \right)(x) - \left(\textup{Id}-\mc{T} \right)(y) \right\|^2 \, .
\end{multline} 
\hfill $\square$
\end{itemize}
\end{definition}

\smallskip
\begin{definition}
A mapping $\mc{T}: \R^n \rightarrow {\R}^n$ is:
\begin{itemize}
\item \textit{Contractive} ({\rm CON}) if there exist $\ell \in [0,1)$ and a norm $\left\| \cdot \right\|$ such that it is an $\ell$-CON in $\mathcal{H}\left( \left\| \cdot \right\| \right)$;
\item \textit{Averaged} ({\rm AVG}) if there exist $\eta \in (0,1)$ and a norm $\left\| \cdot \right\|$ such that it is $\eta$-AVG in $\mathcal{H}\left( \left\| \cdot \right\| \right)$;
\item \textit{strict Pseudo-Contractive} (sPC) if there exists $\kappa \in (0,1)$ and a norm $\left\| \cdot \right\|$ such that it is $\kappa$-sPC in $\mathcal{H}\left( \left\| \cdot \right\| \right)$.

{\hfill $\square$}
\end{itemize}
\end{definition}

%-----------------------------------------------------
%\newpage
\section{Main results: \\ Fixed-point iterations on linear mappings} 
\label{sec:MainRes}

In this section, we provide necessary and sufficient conditions for the convergence of some well-known fixed-point iterations applied on linear operators, i.e.,
%In the remainder of this paper we consider a generic linear mapping
\begin{equation} \label{eq:LinT}
\mathcal{A}:x \mapsto A x, \quad  \text{with } A \in \R^{n\times n}.
\end{equation}

\smallskip

First, we consider the \textit{Banach--Picard} iteration \cite[(1.69)]{bauschke:combettes} on a generic mapping $\cT: \R^n \rightarrow \R^n$, i.e., for all $k \in \N$,
\begin{equation} \label{eq:P-B}
x(k+1) = \cT\left( x(k) \right),
\end{equation}
whose convergence is guaranteed if $\cT$ is averaged, see \cite[Prop. 5.16]{bauschke:combettes}. The next statement shows that averagedness is also a necessary condition when the mapping $\cT$ is linear.

\smallskip
\begin{proposition}[Banach--Picard iteration]
\label{cor:PB-lin}
The following statements are equivalent:
\begin{enumerate}[(i)]
\item $\cA$ in \eqref{eq:LinT} is averaged;
\item the solution to the system
\begin{equation}
x(k+1) = A x(k)
\end{equation}
converges to some $\overline{x} \in \fix (\cA)= \ker (I-A)$.
{\hfill $\square$}
\end{enumerate}
\end{proposition} 
\smallskip

If the mapping $\cT$ is merely nonexpansive, then the sequence generated by the Banach--Picard iteration in \eqref{eq:P-B} may fail to produce a fixed point of $\cT$. For instance, this is the case for $\cT = -\textrm{Id}$. 
In these cases, a relaxed iteration can be used, e.g. the Krasnoselskij--Mann iteration \cite[Equ. (5.15)]{bauschke:combettes}.
Specifically, let us distinguish the case with time-invariant step sizes, known as Krasnoselskij iteration \cite[Chap. 3]{berinde}, and the case with time-varying, vanishing step sizes, known as Mann iteration \cite[Chap. 4]{berinde}. The former is defined by
\begin{equation}\label{eq:Krasno}
x(k+1) = (1-\alpha) x(k) + \alpha \cT\left( x(k) \right) \,,
\end{equation}
for all $k \in \N$, where $\alpha \in (0,1)$ is a constant step size.

The convergence of the discrete-time system in \eqref{eq:Krasno} to a fixed point of the mapping $\cT$ is guaranteed, for any arbitrary $\alpha \in (0,1)$, if $\cT$ is nonexpansive  \cite[Th. 5.15]{bauschke:combettes}, or if $\cT$, defined from a compact, convex set to itself, is strictly pseudo-contractive and $\alpha>0$ is sufficiently small \cite[Theorem 3.5]{berinde}. In the next statement, we show that if the mapping $\cT: \R^n \rightarrow \R^n$ is linear, and $\alpha$ is chosen small enough, then strict pseudo-contractiveness is necessary and sufficient for convergence.

\smallskip
\begin{theorem}[Krasnoselskij iteration]
\label{th:Krasno}
Let $\kappa \in (0,1)$ and $\alpha \in (0,1-\kappa)$.
The following statements are equivalent:
\begin{enumerate}[(i)]
\item $\cA$ in \eqref{eq:LinT} is $\kappa$-strictly pseudo-contractive;
\item the solution to the system
\begin{equation} \label{eq:Krasno-linear}
x(k+1) = (1-\alpha) x(k) + \alpha A x(k)
\end{equation}
converges to some $\overline{x} \in \fix (\cA) = \ker (I-A)$.
%\item the system in \eqref{eq:Krasno-linear} admits a weak polyhedral Lyapunov function.
{\hfill $\square$}
\end{enumerate}
\end{theorem}
\smallskip

In Theorem \ref{th:Krasno}, the admissible step sizes for the Krasnoselskij iteration depend on the parameter $\kappa$ that quantifies the strict pseudo-contractiveness of the mapping $\mathcal{A} = A \, \cdot$. When the parameter $\kappa$ is unknown, or hard to quantify, one can adopt time-varying step sizes, e.g. the Mann iteration:
\begin{equation} \label{eq:MannG}
x(k+1) =  (1-\alpha_k)x(k) + \alpha_k \cT \left( x(k) \right) \,, 
\end{equation}
for all $k \in \N$, where the step sizes $(\alpha_k)_{k \in \N}$ shall be chosen as follows.

\smallskip
\begin{assumption}[Mann sequence]
\label{ass:Mann-sequence}
The sequence $(\alpha_k)_{k \in \N}$ is such that $0 < \alpha_k \leq  \alpha^{\max} < \infty$ for all $k \in \N$, for some $\alpha^{\max}$, $\lim_{k \rightarrow \infty} \alpha_k = 0$ and $\sum_{k=0}^{\infty} \alpha_k = \infty$. 
{\hfill $\square$}
\end{assumption}
\smallskip

The convergence of \eqref{eq:MannG} to a fixed point of the mapping $\cT$ is guaranteed if $\cT$, defined from a compact, convex set to itself, is strictly pseudo-contractive \cite[Theorem 3.5]{berinde}. In the next statement, we show that if the mapping $\cT: \R^n \rightarrow \R^n$ is linear, then strict pseudo-contractiveness is necessary and sufficient for convergence.

\smallskip
\begin{theorem}[Mann iteration] \label{th:Mann}
Let $\left( \alpha_k\right)_{k \in \N}$ be a Mann sequence as in Assumption \ref{ass:Mann-sequence}. The following statements are equivalent:
\begin{enumerate}[(i)]
\item $\cA$ in \eqref{eq:LinT} is strictly pseudocontractive;
\item the solution to 
\begin{equation} \label{eq:Mann-linear}
x(k+1) =  (1-\alpha_k) x(k) + \alpha_k A  x(k)
\end{equation} 
converges to some $\overline{x} \in \fix (\cA)= \ker (I-A)$.
{\hfill $\square$}
\end{enumerate}
\end{theorem}
\smallskip

%----------------------------------------------------------------
%\newpage
\section{Operator-theoretic characterization of linear mappings} \label{sec:OTPLM}

In this section, we characterize the operator-theoretic properties of linear mappings via necessary and sufficient linear matrix inequalities and conditions on the spectrum of the corresponding matrices. We exploit these technical results in Section \ref{sec:main-proofs}, to prove convergence of the fixed-point iterations presented in Section \ref{sec:MainRes}.

\smallskip
\begin{lemma}[Lipschitz continuous linear mapping] \label{lemma:affine-LIP}
Let $\ell > 0$ and $P \in \bbS^n_{\succ 0}$. The following statements are equivalent:
\begin{enumerate}[(i)]
\item $\cA$ in \eqref{eq:LinT} is $\ell$-Lipschitz continuous in $\mathcal{H}\left( \left\| \cdot \right\|_P \right)$;
\item $A^\top P A \preccurlyeq \ell^2 P$.
{\hfill $\square$}
\end{enumerate}
\end{lemma}
\smallskip

\begin{proof}
It directly follows from Definition \ref{def:Lipschitz}.
\end{proof}
\smallskip

\smallskip
\begin{lemma}[Linear contractive/nonexpansive mapping] \label{lem:NE-affine}
Let $\ell \in (0,1)$. The following statements are equivalent:
\begin{enumerate}[(i)]
\item $\cA$ in \eqref{eq:LinT} is an $\ell$-contraction;
\item $\exists P \in \bbS^n_{\succ 0}$ such that $A^\top P A \preccurlyeq \ell^2 P$;

\item the spectrum of $A$ is such that
\begin{equation}
\begin{cases}
\Lambda(A) \subset \ell \, \D_1 \\
\forall \lambda \in \Lambda(A) \cap \bdr( \ell \, \D_1), \ \ \lambda \text{ semi-simple}
\end{cases}
\end{equation} 
\end{enumerate}
If $\ell = 1$, the previous equivalent statements hold if and only if $\cA$ in \eqref{eq:LinT} is nonexpansive.
{\hfill $\square$}
\end{lemma}
\smallskip

\begin{proof} The equivalence between (i) and (ii) follows from Lemma \ref{lemma:affine-LIP}. By the Lyapunov theorem, (iii) holds if and only if the discrete-time linear system $x(k+1) = \frac{1}{\ell} A \,x(k)$ is (at least marginally) stable, i.e., $\Lambda(A) \subset \,\ell \, \D_1$ and the eigenvalues of $A$ on the boundary of the disk, $\Lambda(A) \cap \bdr(\ell \,\D_1)$, are semi-simple. The last statement follows by noticing that an $1$-contractive mapping is nonexpansive.
\end{proof}

\smallskip
\begin{lemma}[Linear averaged mapping] \label{lem:AVG-affine}
Let $\eta \in (0,1)$. The following statements are equivalent:
\begin{enumerate}[(i)]
\item $\cA$ in \eqref{eq:LinT} is $\eta$-averaged;
\item $\exists P \in \bbS^n_{\succ 0}$ such that
\begin{equation*}
A^\top P A \preccurlyeq  \left( 2 \eta - 1\right) P + (1-\eta) \left( A^\top P + P A \right);
\end{equation*}
\item
$\mathcal{A}_{\eta} := A_\eta \, \cdot := \left( 1-\frac{1}{\eta} \right) I\cdot + \frac{1}{\eta}A \cdot$ is nonexpansive;

\item  the spectrum of $A$ is such that
\begin{equation} \label{eq:AVG-eig-cond}
\begin{cases}
\Lambda(A) \subset \D_\eta \\
\forall \lambda \in \Lambda(A) \cap  \bdr (\D_\eta), \ \ \lambda \text{ semi-simple}.
\end{cases} 
\end{equation}
\end{enumerate}
{\hfill $\square$}
\end{lemma}
\smallskip

\begin{proof}
The equivalence (i) $\Leftrightarrow$ (ii) follows directly by inequality \eqref{eq:AVG-ineq} in Definition \ref{def:NE}. By \cite[Prop. 4.35]{bauschke:combettes}, $\cA$ is $\eta$-AVG if and only if the linear mapping $\cA_\eta$ is NE, which proves (i) $\Leftrightarrow$ (iii). To conclude, we show that (iii) $ \Leftrightarrow $ (iv). By Lemma \ref{lem:NE-affine}, the linear mapping $\cA_\eta$ is NE if and only if
\begin{align} \label{eq:AVG-affine-proof-1}
&\begin{cases}
\Lambda(A_\eta) \subset \D_1 \\
\forall \lambda \in \Lambda(A_\eta) \cap  \bdr (\D_1), \ \ \lambda \text{ semi-simple}
\end{cases} \\
\Leftrightarrow &
\begin{cases} \label{eq:AVG-affine-proof-2}
\Lambda(A) \subset (1-\eta)\{1\} + \eta \D_1 = \mathbb{D}_\eta \\
\forall \lambda \in \Lambda(A) \cap  \bdr (\mathbb{D}_\eta), \ \  \lambda \text{ semi-simple}
\end{cases} \
\end{align}
where the equivalence \eqref{eq:AVG-affine-proof-1} $\Leftrightarrow$ \eqref{eq:AVG-affine-proof-2} holds because $\Lambda(A_{\eta}) = (1-\tfrac{1}{\eta})\{1\} + \eta \Lambda(A)$, and because the linear combination with the identity matrix does not alter the geometric multiplicity of the eigenvalues.
\end{proof}

\smallskip
\begin{lemma}[Linear strict pseudocontractive mapping] \label{lem:SPC-affine}
Let $\kappa, \eta \in (0,1)$. The following statements are equivalent:
\begin{enumerate}[(i)]
\item $\cA$ in \eqref{eq:LinT} is $\kappa$-strictly pseudocontractive;
\item $\exists P \in \bbS^n_{\succ 0}$ such that
\begin{equation} \label{eq:LMI-SPC}
(1-\kappa )A^\top P A \preccurlyeq (1+\kappa ) P - \kappa  (A^\top P + PA);
\end{equation}
\item 
$\mathcal{A}_{\kappa}^{\textup{s}} :=  A_\kappa^{\textup{s}} \, \cdot := \kappa I \cdot + (1-\kappa ) A \cdot$ is nonexpansive;

\item the spectrum of $A$ is such that
\begin{equation} \label{eq:SPC-affine-eig}
\begin{cases}
\Lambda(A) \subset \D_{\frac{1}{1-\kappa}} \\
\forall \lambda \in \Lambda(A) \cap  \bdr \left( \mathbb{D}_{\frac{1}{1-\kappa}} \right), \ \  \lambda \text{ semi-simple}
\end{cases}
\end{equation}
\item $\mathcal{A}_{\alpha} := A_\alpha \cdot := (1-\alpha) I \cdot + \alpha A \cdot $ is $\eta$-averaged, with $\alpha = \eta (1-\kappa) \in (0,1)$.
{\hfill $\square$}
\end{enumerate}
\end{lemma}
\smallskip

\begin{proof}
The equivalence (i) $\Leftrightarrow$ (ii) follows directly by inequality \eqref{eq:SPC-ineq} in Definition \eqref{eq:LinT}.
To prove that (ii) $\Leftrightarrow$ (iii), we note that the LMI in \eqref{eq:LMI-SPC} can be recast as
\begin{equation}
\left( \kappa I + (1-\kappa ) A \right)^\top P \left( \kappa I + (1-\kappa ) A \right) \preccurlyeq P ,
\end{equation}
which, by Lemma \ref{lem:NE-affine}, holds true if and only if the mapping $\cA_\kappa^{\textup{s}}$ is NE.\\
(iii) $\Leftrightarrow $ (iv): By Lemma \ref{lem:NE-affine}, $\cA_\kappa^{\textup{s}}$ is NE if and only if
\begin{align} \label{eq:SPC-affine-proof-1}
&\begin{cases}
\Lambda(A_\kappa^{\textup{s}} ) \subset \D_1 \\
\forall \lambda \in \Lambda(A_\kappa^{\textup{s}} ) \cap  \bdr (\D_1), \quad \lambda \text{ semi-simple}
\end{cases}\\
\Leftrightarrow & \label{eq:SPC-affine-proof-2}
\begin{cases} 
\Lambda(A) \subset \left\{ -\frac{\kappa}{1-\kappa } \right\} + \frac{1}{1-\kappa } \, \D_1 =  \D_{\frac{1}{1-\kappa}} \\
\forall \lambda \in \Lambda(A) \cap  \bdr \left( \D_{\frac{1}{1-\kappa}} \right), \ \  \lambda \text{ semi-simple}
\end{cases}
\end{align}
where the equivalence \eqref{eq:SPC-affine-proof-1} $\Leftrightarrow$ \eqref{eq:SPC-affine-proof-2} holds because $\Lambda(A_{\kappa}^{\textup{s}}) =  A_\kappa^{\textup{s}} := \kappa I + (1-\kappa ) A$, and because the linear combination with the identity matrix does not alter the geometric multiplicity of the eigenvalues.
(iii) $\Leftrightarrow$ (v): By Definition \ref{def:NE} and \cite[Prop. 4.35]{bauschke:combettes}, $A_\kappa^{\textup{s}} \cdot $ is NE if and only if $A_\alpha \cdot = (1-\eta)I \cdot + \eta A_\kappa^{\textup{s}} \cdot$ is $\eta$-AVG, for all $\eta \in (0,1)$. Since $\alpha = \eta (1-\kappa)$, $A_\alpha = (1-\eta (1-\kappa)) \textrm{Id} + \eta (1-\kappa) A$, which concludes the proof.
\end{proof}

\section{Proofs of the main results} \label{sec:main-proofs}

\smallskip
\subsection*{Proof of Proposition \ref{cor:PB-lin} (Banach--Picard iteration)}
We recall that, by Lemma \ref{lem:AVG-affine}, $\cA$ is AVG if and only if there exists $\eta \in (0,1)$ such that $\Lambda(A) \subset \mathbb{D}_\eta$ and $\forall \lambda \in \Lambda(A) \cap  \bdr (\mathbb{D}_\eta)$, $\lambda$ is semi-simple and we notice that  $\D_\eta \cap \D_1 = \{ 1 \}$ for all $ \eta \in (0,1)$. Hence $\cA$ is averaged if and only if the eigenvalues of $A$ are strictly contained in the unit circle except for the eigenvalue in $\lambda = 1$ which, if present, is semi-simple. The latter is a necessary and sufficient condition for the convergence of $x(k+1) = A \, x(k)$, by Lemma \ref{lem:conv-eig}.
{\hfill $\blacksquare$}

\smallskip
\subsection*{Proof of Theorem \ref{th:Krasno} (Krasnoselskij iteration) }
(i) $\Leftrightarrow $ (ii): By Lemma \ref{lem:SPC-affine}, $\cA$ is $k$-sPC if and only if $(1-\alpha)\Id + \alpha \cA$ is $\eta$-AVG, with  $\alpha = \eta (1-\kappa)$ and $\eta \in (0,1)$; therefore, if and only if $(1-\alpha)\Id + \alpha \cA$ is AVG with $\alpha \in (0,1-\kappa)$. By proposition \eqref{cor:PB-lin}, the latter is equivalent to the global convergence of the Banach--Picard iteration applied on $(1-\alpha)\Id + \alpha \cA$, which corresponds to the Krasnoselskij iteration on $\cA$, with $\alpha \in (0,1-\kappa)$. 
%(ii) $\Leftrightarrow$ (iii): The proof follows by \cite[Th. 4.50]{blanchini:miani}, since the system in \eqref{eq:Krasno-linear} converges, hence at least marginally stable \cite[Th. 4.50 (i)]{blanchini:miani}, and the matrix $(1-\alpha) I + \alpha A$ has eigenvalues strictly inside the unit disk, or in $1$ \cite[Th. 4.50 (ii, $\theta=0$)]{blanchini:miani}.
{\hfill $\blacksquare$}

\smallskip
\subsection*{Proof of Theorem \ref{th:Mann} (Mann iteration)}
Proof that (i) $\Rightarrow$ (ii): 
Define the bounded sequence ${\beta}_k := \tfrac{1}{\epsilon} \alpha_k > 0$, for some $\epsilon > 0$ to be chosen. Thus, $x(k+1) = (1-\alpha_k) x(k) + \alpha_k A x(k) = (1- \epsilon \beta_k) x(k) + \epsilon \beta_k A x(k) = (1-\beta_k) x(k) + \beta_k \left(  (1-\epsilon) I + \epsilon A \right) x(k) $. Since $A \cdot$ is sPC, we can choose $\epsilon>0$ small enough such that $B := (1-\epsilon) \Id + \epsilon A \cdot$ is NE, specifically, we shall choose $\epsilon < \min\left\{ \tfrac{1}{|\lambda|} \mid  \lambda \in \Lambda(A)\setminus \{ 1 \} \right\}$. Note that $0 \in \textrm{fix}(A) = \textrm{fix}(B) \neq \varnothing$. Since $\infty = \sum_{k=0}^{\infty} \alpha_k = \epsilon \sum_{k=0}^{\infty} \beta_k$, we have that $\lim_{k \rightarrow \infty} \beta_k = 0$, hence $\forall \epsilon > 0$, $\exists \bar{k} \in \N$ such that $\beta_k \leq 1$ for all $k \geq \bar{k}$. Moreover, since $\sum_{k=0}^{\bar{k}} \beta_k < \infty$, for all $x(0) \in \R^n$, we have that the solution $x(\bar{k})$ is finite. Therefore, we can define $h := k - \bar{k} \in \N$ for all $k \geq \bar{k}$, $y(0) := x(\bar{k})$ and $y(h+1) := x( h + \bar{k} + 1 )$ for all $h \geq 0$. The proof then follows by applying \cite[Th. 5.14 (iii)]{bauschke:combettes} to the Mann iteration $y(h+1) = (1-\beta_h) y(h) + \beta_h B y(h)$.

Proof that (ii) $\Rightarrow$ (i): 
For the sake of contradiction, suppose that $A$ is not sPC, i.e., at least one of the following facts must hold: 1) $A$ has an eigenvalue in $1$ that is not semi-simple; 2) $A$ has a real eigenvalue greater than $1$; 3) $A$ has a pair of complex eigenvalues $\sigma \pm j \omega$, with $\sigma \geq 1$ and $\omega > 0$. 
We show next that each of these three facts implies non-convergence of \eqref{eq:MannG}. Without loss of generality (i.e., up to a linear transformation), we can assume that $A$ is in Jordan normal form. 

1) $A$ has an eigenvalue in $1$ that is not semi-simple. Due to (the bottom part of) the associated Jordan block, the vector dynamics in \eqref{eq:MannG} contain the two-dimensional linear time-varying dynamics 
\vspace{-0.15cm}
\begin{align*}
y(k+1) &= 
\left(
(1-\alpha_k) 
\left[ 
\begin{matrix}
1 & 0 \\
0 & 1
\end{matrix}
\right] + \alpha_k 
\left[ 
\begin{matrix}
1 & 1 \\
0 & 1
\end{matrix}
\right] \right) y(k) \\
&= \left[ 
\begin{matrix}
1 & \alpha_k \\
0 & 1
\end{matrix}
\right] y(k).
\end{align*}
For $y_2(0) := c > 0$, we have that the solution $y(k)$ is such that $y_2(k) = y_2(0) > 0 $ and $y_1(k+1) = y_1(k) + c$, which implies that $y_1(k) = y_1(0) + k \, c$. Thus, $x(k)$ diverges and we have a contradiction.

2) Let $A$ has a real eigenvalue equal to $1 + \epsilon > 1$. Again due to (the bottom part of) the associated Jordan block, the vector dynamics in \eqref{eq:MannG} must contain the scalar dynamics $s(k+1) = (1-\alpha_k) s(k) + \alpha_k ( 1 + \epsilon ) s(k) = (1 + \epsilon \,  \alpha_k ) s(k)$. The solution then reads as $s(k+1) = \left( \prod_{h=0}^{k} ( 1 + \epsilon \, \alpha_h ) \right) s(0) $. Now, since $\epsilon \, \alpha_h > 0$, it holds that $\prod_{h=0}^{k} ( 1 + \epsilon \, \alpha_h ) \geq \epsilon \sum_{h=0}^{k} \alpha_h  = \infty$, by Assumption \ref{ass:Mann-sequence}. Therefore, $s(k)$ and hence $x(k)$ diverge, and we reach a contradiction.

3) $A$ has a pair of complex eigenvalues $\sigma \pm j \omega$, with $\sigma = 1+\epsilon \geq 1$ and $\omega > 0$. Due to the structure of the associated Jordan block, the vector dynamics in \eqref{eq:MannG} contain the two-dimensional dynamics 
\vspace{-0.15cm}
\begin{align*}
z(k+1) & = 
\left(
(1-\alpha_k) 
\left[ 
\begin{matrix}
1 & 0 \\
0 & 1
\end{matrix}
\right] + \alpha_k 
\left[ 
\begin{matrix}
\sigma & -\omega \\
\omega & \phantom{-}\sigma
\end{matrix}
\right] \right) z(k) \\
&= 
\left[ 
\begin{matrix}
1+\epsilon \, \alpha_k & -\omega \alpha_k \\
\omega \alpha_k & 1+\epsilon \, \alpha_k
\end{matrix}
\right] z(k).
\end{align*}

Now, we define $\rho_k := \sqrt{(1+\epsilon \alpha_k)^2 + \omega^2 \alpha_k^2 } \geq \sqrt{1 + \omega^2 \alpha_k^2} > 1$, and the angle $\theta_k > 0$ such that $\cos(\theta_k) = (1+\epsilon \alpha_k)/\rho_k$ and $\sin(\theta_k) = (\omega \alpha_k)/\rho_k$, i.e., $\theta_k = \textrm{atan}\left( \frac{\omega \alpha_k}{ 1 + \epsilon \alpha_k}\right)$. Then, we have that $ z(k+1) = \rho_k R(\theta_k) z(k)$,
%where $R(\cdot) := \left[ 
%\begin{smallmatrix}
%\cos(\cdot) & -\sin(\cdot) \\
%\sin(\cdot) & \phantom{-}\cos(\cdot)
%\end{smallmatrix}
%\right]$ is the rotation operator. 
hence, the solution $z(k)$ reads as 
$$\textstyle z(k+1) = \left( \prod_{h=0}^{k} \rho_h \right) R\left( \sum_{h=0}^{k} \theta_h \right) z(0).$$

Since $\left\| R(\cdot)\right\| = 1$, if the product $\left( \prod_{h=0}^{\infty} \rho_h \right)$ diverges, then $z(k)$ and hence $x(k)$ diverge as well. Thus, let us assume that the product $\left( \prod_{h=0}^{\infty} \rho_h \right)$ converges. By the limit comparison test, the series $\sum_{h=0}^{\infty} \theta_h =  \sum_{h=0}^{\infty} \textrm{atan}\left( \frac{\omega \alpha_h}{ 1 + \epsilon \alpha_h} \right)$ converges (diverges) if and only the series $\sum_{h=0}^{\infty} \frac{\omega \alpha_h}{ 1 + \epsilon \alpha_h} $ converges (diverges). The latter diverges since $ \sum_{h=0}^{\infty} \frac{\omega \alpha_h}{ 1 + \epsilon \alpha_h} \geq \omega \sum_{h=0}^{\infty} \frac{\alpha_h}{ 1 + \epsilon \alpha^{\max}} = \frac{\omega}{1 + \epsilon \alpha^{\max}} \sum_{h=0}^{\infty} \alpha_h = \infty$. It follows that $\sum_{h=0}^{\infty} \theta_h$ diverges, hence $z(k)$ keeps rotating indefinitely, which is a contradiction.
{\hfill $\blacksquare$}
\smallskip

%----------------------------------------------------------------

%-------------------------------------------------------
\section{Application to multi-agent linear systems} \label{sec:App}

\subsection{Consensus via time-varying Laplacian dynamics}
\label{sec:App-Cons}
We consider a connected graph of $N$ nodes, associated with $N$ agents seeking consensus, with Laplacian matrix $L \in \R^{N \times N}$. To solve the consensus problem, we study the following discrete-time linear time-varying dynamics:
\begin{subequations}\label{eq:tv-consensus}
\begin{align}
\bs{x}(k+1) & = \bs{x}(k) - \alpha_k L \bs{x}(k) \\
            & = (1-\alpha_k)\bs{x}(k) + \alpha_k (I-L) \bs{x}(k) \label{eq:tv-consensus-Mann} \,,
\end{align}
\end{subequations}
where $\bs{x}(k) := \left[ x_1(k), \ldots, x_N(k) \right]^\top \in \R^N$ and, for simplicity, the state of each agent is a scalar variable, $x_i \in \R$.

Since the dynamics in \eqref{eq:tv-consensus} have the structure of a Mann iteration, in view of Theorem \ref{th:Mann}, we have the following result.

\smallskip
\begin{corollary}
Let $\left( \alpha_k\right)_{k \in \N}$ be a Mann sequence.
The system in \eqref{eq:tv-consensus} asymptotically reaches consensus, i.e., the solution $\bs{x}(k)$  to \eqref{eq:tv-consensus} converges to $\overline{x} \, \bs{1}_N$, for some $\overline{x} \in \R$.
{\hfill $\square$}
\end{corollary}
\smallskip

\begin{proof}
Since the graph is connected, $L$ has one (simple) eigenvalue at $0$, and $N-1$ eigenvalues with strictly-positive real part. Therefore, the matrix $I-L$ in \eqref{eq:tv-consensus-Mann} has one simple eigenvalue in $1$ and $N-1$ with real part strictly less than $1$. By Lemma \ref{lem:SPC-affine}, $(I-L)(\cdot)$ is sPC and by Theorem \ref{th:Mann}, $\bs{x}(k)$ globally converges to some $\overline{\bs{x}} \in \textrm{fix}(I-L) = \textrm{zer}(L)$, i.e., $L \overline{\bs{x}} = \bs{0}_N$. Since $L$ is a Laplacian matrix, $L \overline{\bs{x}} = \bs{0}_N$ implies consensus, i.e., $\overline{\bs{x}} = \overline{x} \, \bs{1}_N$, for some $\overline{x} \in \R$.
\end{proof}
\smallskip

We emphasize that via \eqref{eq:tv-consensus}, consensus is reached without assuming that the agents know the algebraic connectivity of the graph, i.e., the strictly-positive Fiedler eigenvalue of $L$. We have only assumed that the agents agree on a sequence of vanishing, bounded, step sizes, $\alpha_k$. However, we envision that agent-dependent step sizes can be used as well, e.g. via matricial Mann iterations, see \cite[\S 4.1]{berinde}.

Let us simulate the time-varying consensus dynamics in \eqref{eq:tv-consensus} for a graph with $N=3$ nodes, adjacency matrix $A = \left[ a_{i,j} \right]$ with $a_{1,2} = a_{1,3} = 1/2$, $a_{2,3} = a_{3,1} = 1$, hence with Laplacian matrix
$$ L = D_{\textup{out}} - A =
\left[
\begin{smallmatrix}
\phantom{-}1 & -1/2 &  -1/2 \\
\phantom{-}0 & 1 & -1\\
-1 & 0 & \phantom{-}1
\end{smallmatrix}
\right].
 $$
 
\begin{figure}
\includegraphics[width=\columnwidth]{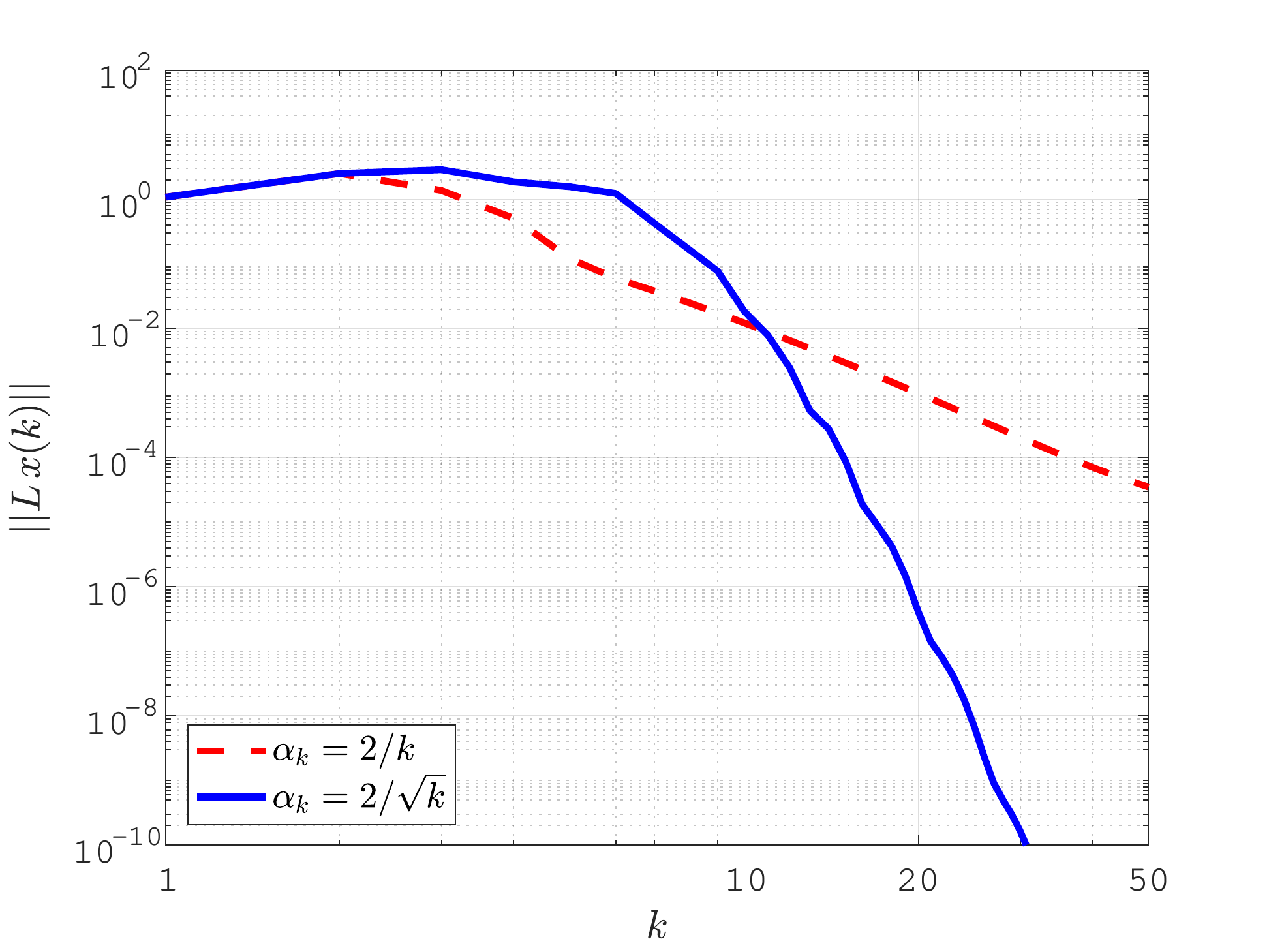}
\caption{Disagreement vector norm versus discrete time. Since the mapping $\Id - L \, \cdot$ is strictly pseudocontractive, consensus is asymptotically reached.}
\label{fig:consensus}
\end{figure}

We note that $L$ has eigenvalues $\Lambda(L)=\left\{ 0 , \tfrac{3}{2} \pm j \tfrac{1}{2} \right\}$. Since we do not assume that the agents known about the connectivity of the graph, we simulate with step sizes that are initially larger than the maximum constant-step value for which convergence would hold. In Fig. \ref{fig:consensus}, we compare the norm of the disagreement vectors, $\left\| L x(k)\right\|$,  obtained with two different Mann sequences, $\alpha_k = 2/k$ and $\alpha_k = 2/\sqrt{k}$, respectively. We observe that convergence with small tolerances is faster in the latter case with larger step sizes.

\subsection{Two-player zero-sum linear-quadratic games: \\ Non-convergence of projected pseudo-gradient dynamics} \label{sec:App-Game}

We consider two-player zero-sum games with linear-quadratic structure, i.e., we consider $N=2$ agents, with cost functions $f_1(x_1, x_2) := x_1^\top C x_2$ and $f_2(x_1, x_2) := -x_2^\top C^\top x_1$, respectively, for some square matrix $C = C^\top \neq 0$. In particular, we study discrete-time dynamics for solving the Nash equilibrium problem, that is the problem to find a pair $\left( x_1^*, x_2^* \right)$ such that:
\begin{equation*}
\left\{
\begin{array}{rl}
x_1^* & \in \, \underset{ x_1 \in \R^n }{\textrm{argmin}} f_1( x_1, x_2^* ) \\
x_2^* & \in \, \underset{ x_2 \in \R^n }{\textrm{argmin}} f_2( x_1^*, x_2 ).
\end{array}
\right.
\end{equation*}

A classic solution approach is the pseudo-gradient method, namely the discrete-time dynamics 
\begin{subequations}\label{eq:Mann-two-player-games}
\begin{align} 
\bs{x}(k+1) & =  \bs{x}(k) -\alpha_k F \bs{x}(k)\\
            & = ( 1 - \alpha_k ) \bs{x}(k) + \alpha_k (I-F) \bs{x}(k) \,,
\end{align}
\end{subequations}
where $F \, \cdot $ is the so-called pseudo-gradient mapping of the game, which in our case is defined as
$$ \mathcal{F}( \bs{x} ) := \left[ 
\begin{matrix}
\nabla_{x_1} f_1(x_1, x_2) \\
\nabla_{x_2} f_2(x_1, x_2)
\end{matrix}
\right] = \left[ 
\begin{matrix}
\phantom{-}C x_2 \\
-C x_1
\end{matrix}
\right] = 
\underbrace{\left[ 
\begin{matrix}
\phantom{-}0 & 1 \\
-1 & 0
\end{matrix}
\right] \otimes C}_{ =: F } \, \bs{x} \,,$$
and $(\alpha_k)_{k \in \N}$ is a sequence of vanishing step sizes, e.g. a Mann sequence.
In our case, $\left( x_1^{*}, x_2^{*} \right)$ is a Nash equilibrium  if and only if $ \left[ x_1^{*}  \, ; x_2^{*} \right] \in \textrm{fix}\left( \Id - \mathcal{F}\right) = \zer\left( \mathcal{F}\right)$ \cite[Th. 1]{belgioioso:grammatico:17cdc}.

By Theorem \ref{th:Mann}, convergence of the system in \eqref{eq:Mann-two-player-games} holds if and only if $I-F$ is strictly pseudocontractive. In the next statement, we show that this is not the case for $F$ in \eqref{eq:Mann-two-player-games}.

\smallskip
\begin{corollary}
Let $(\alpha_k)_{k \in \N}$ be a Mann sequence and $C=C^\top \neq 0$. The system in \eqref{eq:Mann-two-player-games} does not globally converge.
{\hfill $\square$}
\end{corollary}

\smallskip
\begin{proof}
It follows by Lemma \ref{lem:SPC-affine} that the mapping $\Id - F \cdot$ is strictly pseudocontractive if and only if the eigenvalues of $F$ either have strictly-positive real part, or are semi-simple and equal to $0$. 
Since $\Lambda\left( \left[ \begin{smallmatrix} \phantom{-}0 & 1 \\ -1 & 0 \end{smallmatrix} \right] \right) = \{ \pm j \}$, we have that the eigenvalues of $F = \left[ \begin{smallmatrix} \phantom{-}0 & 1 \\ -1 & 0 \end{smallmatrix} \right] \otimes C$ are either with both positive and negative real part, or on the imaginary axis and not equal to $0$, or equal to $0$ are not semi-simple. Therefore, $\Id - F \cdot$ is not strictly pseudocontractive and the proof follows by Theorem \ref{th:Mann}.
\end{proof}
\smallskip

\begin{figure}[!h]
\includegraphics[width=\columnwidth]{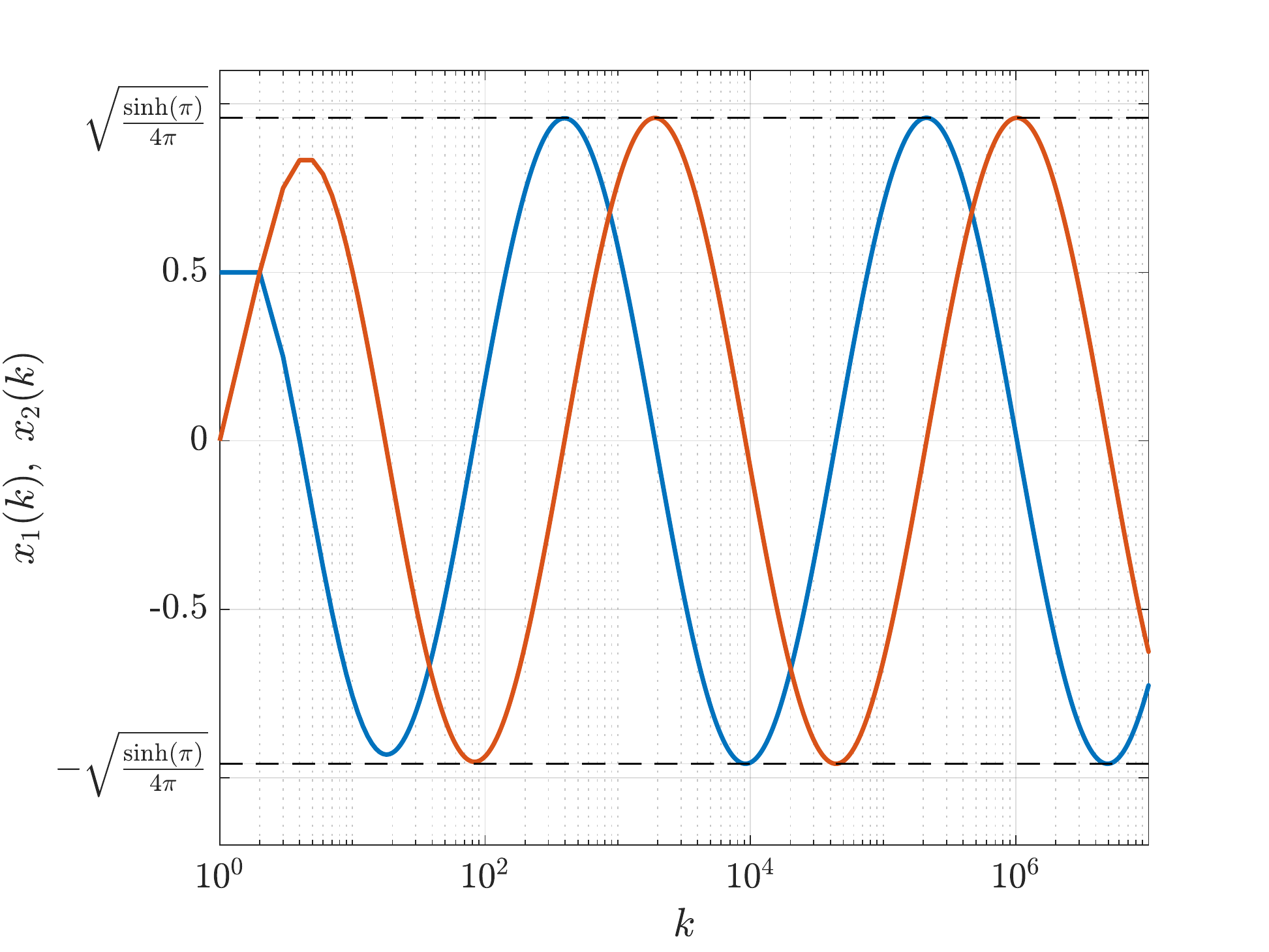}
\caption{Solution to the discrete-time system in \eqref{eq:Mann-two-player-games} in semi-log scale. 
The lack of strict pseudo-contractiveness causes persistent oscillations.}
\label{fig:exampleB}
\end{figure}

Let us numerically simulate the discrete-time system in \eqref{eq:Mann-two-player-games}, with the following parameters: $n=1$, $C = 1$, $x_1(0) = 1/2$, $x_2(0)=0$, and $\alpha_k = 1/(k+1)$ for all $k \in \N$. Figure \ref{fig:exampleB} shows persistent oscillations, due to the lack of strict pseudo-contractiveness of $I-F$. In fact, $\Lambda(I-F) = \{ 1 \pm j \}$. The example provides a non-convergence result: pseudo-gradient methods do not ensure global convergence in convex games with (non-strictly) monotone pseudo-gradient mapping, not even with vanishing step sizes and linear-quadratic structure. 
%A similar conclusion can be derived for monotone games with local and coupling constraints.

%----------------------------------------------------------------

\section{Conclusion and outlook} \label{sec:Concl}

Convergence in discrete-time linear systems can be equivalently characterized via operator-theoretic notions. 
Remarkably, the time-varying Mann iteration applied on linear mappings converges if and only if the considered linear operator is strictly pseudocontractive. This result implies that Laplacian-driven linear time-varying consensus dynamics with Mann step sizes do converge. It also implies that projected pseudo-gradient dynamics for Nash equilibrium seeking in monotone games do not necessarily converge. 

Future research will focus on studying convergence of other, more general, linear fixed-point iterations and of discrete-time linear systems with uncertainty, e.g. polytopic.

%\balance
\bibliographystyle{IEEEtran}
\bibliography{library}

\begin{thebibliography}{10}
\providecommand{\url}[1]{#1}
\csname url@samestyle\endcsname
\providecommand{\newblock}{\relax}
\providecommand{\bibinfo}[2]{#2}
\providecommand{\BIBentrySTDinterwordspacing}{\spaceskip=0pt\relax}
\providecommand{\BIBentryALTinterwordstretchfactor}{4}
\providecommand{\BIBentryALTinterwordspacing}{\spaceskip=\fontdimen2\font plus
\BIBentryALTinterwordstretchfactor\fontdimen3\font minus
  \fontdimen4\font\relax}
\providecommand{\BIBforeignlanguage}[2]{{%
\expandafter\ifx\csname l@#1\endcsname\relax
\typeout{** WARNING: IEEEtran.bst: No hyphenation pattern has been}%
\typeout{** loaded for the language `#1'. Using the pattern for}%
\typeout{** the default language instead.}%
\else
\language=\csname l@#1\endcsname
\fi
#2}}
\providecommand{\BIBdecl}{\relax}
\BIBdecl

\bibitem{olfati-saber:fax:murray:10}
R.~Olfati-Saber, A.~Fax, and R.~Murray, ``Consensus and cooperation in
  networked multi-agent systems,'' \emph{Proc. of the IEEE}, vol.~95, no.~1,
  pp. 215--233, 2010.

\bibitem{nedic:ozdaglar:parrillo:10}
A.~Nedi\'c, A.~Ozdaglar, and P.~Parrillo, ``Constrained consensus and
  optimization in multi-agent networks,'' \emph{IEEE Trans. on Automatic
  Control}, vol.~55, no.~4, pp. 922--938, 2010.

\bibitem{yi:pavel:17cdc}
P.~Yi and L.~Pavel, ``A distributed primal-dual algorithm for computation of
  generalized {Nash} equilibria via operator splitting methods,'' \emph{Proc.
  of the IEEE Conf. on Decision and Control}, pp. 3841--3846, 2017.

\bibitem{dorfer:simpson-porco:bullo:16}
F.~D\"{o}rfler, J.~Simpson-Porco, and F.~Bullo, ``Breaking the hierarchy:
  Distributed control and economic optimality in microgrids,'' \emph{IEEE
  Trans. on Control of Network Systems}, vol.~3, no.~3, pp. 241--253, 2016.

\bibitem{dorfler:grammatico:17}
F.~D\"{o}rfler and S.~Grammatico, ``Gather-and-broadcast frequency control in
  power systems,'' \emph{Automatica}, vol.~79, pp. 296--305, 2017.

\bibitem{mohsenian-rad:10}
A.-H. Mohsenian-Rad, V.~Wong, J.~Jatskevich, R.~Schober, and A.~Leon-Garcia,
  ``Autonomous demand-side management based on game-theoretic energy
  consumption scheduling for the future smart grid,'' \emph{IEEE Trans. on
  Smart Grid}, vol.~1, no.~3, pp. 320--331, 2010.

\bibitem{jaina:walrand:10}
R.~Jaina and J.~Walrand, ``An efficient {Nash}-implementation mechanism for
  network resource allocation,'' \emph{Automatica}, vol.~46, pp. 1276Ð--1283,
  2010.

\bibitem{barrera:garcia:15}
J.~Barrera and A.~Garcia, ``Dynamic incentives for congestion control,''
  \emph{IEEE Trans. on Automatic Control}, vol.~60, no.~2, pp. 299--310, 2015.

\bibitem{Ghaderi2014}
J.~Ghaderi and R.~Srikant, ``Opinion dynamics in social networks with stubborn
  agents: Equilibrium and convergence rate,'' \emph{Automatica}, vol.~50, pp.
  3209--Ð3215, 2014.

\bibitem{etesami:basar:15}
S.~R. Etesami and T.~Ba\c{s}ar, ``Game-theoretic analysis of the
  hegselmann-krause model for opinion dynamics in finite dimensions,''
  \emph{IEEE Trans. on Automatic Control}, vol.~60, no.~7, pp. 1886--Ð1897,
  2015.

\bibitem{martinez:bullo:cortes:frazzoli:07}
S.~Mart\'{i}nez, F.~Bullo, J.~Cort\'{e}s, and E.~Frazzoli, ``On synchronous
  robotic networks -- {Part} i: Models, tasks, and complexity,'' \emph{IEEE
  Trans. on Automatic Control}, vol.~52, pp. 2199--2213, 2007.

\bibitem{sankovic:johansson:stipanovic:12}
M.~Stankovi\'{c}, K.~Johansson, and D.~Stipanovi\'{c}, ``Distributed seeking of
  {Nash} equilibria with applications to mobile sensor networks,'' \emph{IEEE
  Trans. on Automatic Control}, vol.~57, no.~4, pp. 904--919, 2012.

\bibitem{berinde}
V.~Berinde, \emph{Iterative Approximation of Fixed Points}.\hskip 1em plus
  0.5em minus 0.4em\relax Springer, 2007.

\bibitem{bauschke:combettes}
H.~H. Bauschke and P.~L. Combettes, \emph{Convex analysis and monotone operator
  theory in {Hilbert} spaces}.\hskip 1em plus 0.5em minus 0.4em\relax Springer,
  2010.

\bibitem{ryu:boyd:16}
E.~K. Ryu and S.~Boyd, ``A primer on monotone operator methods,'' \emph{Appl.
  Comput. Math.}, vol.~15, no.~1, pp. 3--43, 2016.

\bibitem{grammatico:parise:colombino:lygeros:16}
S.~Grammatico, F.~Parise, M.~Colombino, and J.~Lygeros, ``Decentralized
  convergence to {Nash} equilibria in constrained deterministic mean field
  control,'' \emph{IEEE Trans. on Automatic Control}, vol.~61, no.~11, pp.
  3315--3329, 2016.

\bibitem{grammatico:17}
S.~Grammatico, ``Dynamic control of agents playing aggregative games with
  coupling constraints,'' \emph{IEEE Trans. on Automatic Control}, vol.~62,
  no.~9, pp. 4537 -- 4548, 2017.

\bibitem{giselsson:boyd:17}
P.~Giselsson and S.~Boyd, ``Linear convergence and metric selection for
  {Douglas--Rachford} splitting and {ADMM},'' \emph{IEEE Transactions on
  Automatic Control}, vol.~62, no.~2, pp. 532--544, 2017.

\bibitem{belgioioso:grammatico:17cdc}
G.~Belgioioso and S.~Grammatico, ``Semi-decentralized {Nash} equilibrium
  seeking in aggregative games with coupling constraints and non-differentiable
  cost functions,'' \emph{IEEE Control Systems Letters}, vol.~1, no.~2, pp.
  400--405, 2017.

\bibitem{grammatico:17tcns}
S.~Grammatico, ``Proximal dynamics in multi-agent network games,'' \emph{IEEE
  Trans. on Control of Network Systems,
  \textup{\texttt{\url{https://doi.org/10.1109/TCNS.2017.2754358}}}}, 2018.

\end{thebibliography}

\end{document}